\documentclass[12pt]{amsart}

\usepackage{amscd, amssymb}

\input{prepictex}
\input{pictex}
\input{postpictex}

\newtheorem{thm}{Theorem}[section]
\newtheorem{cor}[thm]{Corollary}
\newtheorem{lem}[thm]{Lemma}
\newtheorem{prop}[thm]{Proposition}

\theoremstyle{definition}
\newtheorem{dfn}[thm]{Definition}

\theoremstyle{remark}
\newtheorem{rmk}[thm]{Remark}
\newtheorem{rmks}[thm]{Remarks}

\newtheorem{examples}[thm]{Examples}

\numberwithin{equation}{section}

\begin{document}

\title{Flow equivalence of graph algebras}

\author{Teresa Bates and David Pask}

\address{University of New South Wales, Sydney NSW 2052,
AUSTRALIA} \email{teresa@maths.unsw.edu.au}

\address{The University of Newcastle, Callaghan, NSW 2308}
\email{davidp@maths.newcastle.edu.au}

\date{December 2002}

\begin{abstract}
This paper explores the effect of various graphical constructions
upon the associated graph $C^*$-algebras. The graphical
constructions in question arise naturally in the study of flow
equivalence for topological Markov chains.
 We prove that out-splittings  give
rise to isomorphic graph algebras, and in-splittings give rise to
strongly Morita equivalent $C^*$-algebras.  We generalise the
notion of a delay as defined in \cite{d} to form in-delays and
out-delays. We prove that these constructions give rise to Morita
equivalent graph $C^*$-algebras.  We provide examples which
suggest that our results are the most general possible in the
setting of the $C^*$-algebras of arbitrary directed graphs.
\end{abstract}

\maketitle

\section{Introduction}

The purpose of this paper is to describe various constructions on
a directed graph which give rise to equivalences between the
associated graph $C^*$-algebras. The graphical constructions in
question all have their roots in the theory of flow equivalence
for topological Markov chains. Our results will unify the work of
several authors over the last few years who have studied similar
constructions for Cuntz-Krieger algebras and, more recently, graph
$C^*$-algebras (see \cite{ck,mrs,ash,d,ds,bb,b1} amongst others).

The motivation for the graphical constructions we use lies in the
theory of subshifts of finite type. A shift space $( \textsf{X} ,
\sigma )$ over a finite alphabet $\mathcal{A}$ is a compact subset
$\textsf{X}$ of $\mathcal{A}^{\bf Z}$ invariant under the shift
map $\sigma$. To a directed graph $E$ with finitely many edges and
no sources or sinks one may associate a shift space
$\textsf{X}_E$, called the edge shift of $E$, whose alphabet is
the edge set of $E$ (see \cite[Definition 2.2.5]{lm}). Edge shifts
are examples of subshifts of finite type. Alternatively, to every
square $0$-$1$ matrix $A$ with no zero rows or columns one may
associate a subshift of finite type $\textsf{X}_A$ (see
\cite[Definition 2.3.7]{lm}).

Two important types of equivalence between shift spaces are
conjugacy and flow equivalence. Shift spaces $( \textsf{X} ,
\sigma_{\textsf{X}} )$ and $( \textsf{Y} , \sigma_{\textsf{Y}} )$
are conjugate ($\textsf{X} \cong \textsf{Y}$) if there is an
isomorphism $\phi : \textsf{X} \to \textsf{Y}$ such that
$\sigma_{\textsf{Y}} \circ \phi = \phi \circ \sigma_{\textsf{X}}$.
The suspension of $( \textsf{X} , \sigma_{\textsf{X}} )$ is
\[
\textsf{SX} := ( \textsf{X} \times {\bf R} ) / [ ( x, t+1 ) \sim (
\sigma_{\textsf{X}} (x) , t ) ] .
\]

\noindent Shift spaces $( \textsf{X} , \sigma_{\textsf{X}} )$ and
$( \textsf{Y} , \sigma_{\textsf{Y}} )$ are flow equivalent
($\textsf{X} \sim \textsf{Y}$) if there is a homeomorphism between
$\textsf{SX}$ and $\textsf{SY}$ preserving the orientation of flow
lines.

By \cite[Proposition 2.3.9]{lm} every subshift of finite type is
conjugate to an edge shift $\textsf{X}_E$ for some  directed graph
$E$. Since the edge connectivity matrix $B_E$ of $E$ is a $0$-$1$
matrix such that $\textsf{X}_E \cong \textsf{X}_{B_E}$, every
subshift of finite type is conjugate to a shift described by a
$0$-$1$ matrix. Conversely, every shift described by a $0$-$1$
matrix $A$ is conjugate to an edge shift: Let $E_A$ be the
directed graph with vertex connectivity matrix $A$, then
$\textsf{X}_{E_A}$ is conjugate to $\textsf{X}_A$ (see
\cite[Exercise 1.5.6, Proposition 2.3.9]{lm}). Hence subshifts of
finite type are edge shifts or shifts associated to $0$-$1$
matrices.

Conjugacy and flow equivalence for subshifts of finite type may be
expressed in terms of $0$-$1$ matrices: An elementary strong shift
equivalence between square $0$-$1$ matrices $A , B$ is a pair
$(R,S)$ of $0$-$1$ matrices such that $A=RS$ and $B=SR$. We say
$A$ and $B$ are strong shift equivalent if there is a chain of
elementary strong shift equivalences from $A$ to $B$. From
\cite[Theorem A]{wi} (see also \cite[Theorem 7.2.7]{lm})
$\textsf{X}_A \cong \textsf{X}_B$ if and only if $A$ and $B$ are
strong shift equivalent. By \cite{ps} $\textsf{X}_A \sim
\textsf{X}_B$ if and only if $A$ and $B$ are related via a chain
of elementary strong shift equivalences and certain matrix
expansions. Both these matrix operations have graphical
interpretations: Following \cite[Theorem 2.4.14 and Exercise
2.4.9]{lm} an elementary strong shift equivalence corresponds to
either an in- or out-splitting of the corresponding graphs.
Following \cite[\S 3]{d}, the matrix expansions in \cite{ps}
correspond to an out-delay of the corresponding graph.

To a $0$-$1$ matrix $A$ with $n$ non-zero rows and columns is
associated a $C^*$-algebra generated by partial isometries $\{ S_i
\}_{i=1}^n$ with mutually orthogonal ranges satisfying
\[
S_j^* S_j = \sum_{i=1}^n A(i,j) S_i S_i^*.
\]

\noindent If $A$ satisfies condition (I) the Cuntz-Krieger algebra
$\mathcal{O}_A$ is unique up to isomorphism. Results about
Cuntz-Krieger algebras may be expressed in terms of the 
directed graph $E_A$ associated to $A$ (see \cite{ew} and
\cite{fw}, for instance). More recent results are expressed
entirely in terms of $E_A$ (see \cite{kprr,kpr} amongst others).

To a row-finite directed graph $E$ with finitely many edges and no
sources or sinks is associated the universal $C^*$-algebra, $C^*
(E)$ generated by partial isometries $\{ s_e : e \in E^1 \}$ with
mutually orthogonal ranges satisfying
\[
s_e^*s_e = \sum_{s(f)=r(e)} s_f s_f^* .
\]

\noindent If $A$ is a square $0$-$1$ matrix which satisfies
condition (I), and $E_A$ is the associated directed graph, then
$\mathcal{O}_A \cong C^* ( E_A )$ (see \cite[Proposition
4.1]{mrs}). On the other hand, if $E$ satisfies condition (L),
then the associated edge connectivity matrix $B_E$ satisfies
condition (I) and $C^* (E) \cong \mathcal{O}_{B_E}$ (see
\cite[Proposition 4.1]{kprr}).  There are similar equivalences
between Cuntz-Krieger algebras associated to infinite $0$-$1$
matrices which are row-finite and certain row-finite directed
graphs (see \cite{pr}, \cite[Theorem 3.1]{bprs}).

By \cite[Proposition 2.17, Theorems 3.8 and 4.1]{ck} if $A$ and
$B$ satisfy condition (I) and $\textsf{X}_A \cong \textsf{X}_B$
then $\mathcal{O}_A \cong \mathcal{O}_B$; moreover, if
$\textsf{X}_A \sim \textsf{X}_B$ then $\mathcal{O}_A$ is stably
isomorphic to $\mathcal{O}_B$. The aim of this paper is to show
that the graphical procedures involved in flow equivalence and
conjugacy for edge shifts may be applied to arbitrary graphs, and
give rise to isomorphisms or Morita equivalences of the
corresponding graph $C^*$-algebras. Initial results in this
direction were proved in \cite{d} for graphs with no sinks and
finitely many vertices: Out-splittings lead to isomorphisms of the
underlying graph $C^*$-algebras whilst the $C^*$-algebra of an
in-split graph is isomorphic to the $C^*$-algebra of a certain
out-delayed graph. Further partial results may be found in
\cite{ds} and \cite{bb}.

 The paper is organised as follows. Section
\ref{basic} describes the $C^*$-algebra of any directed graph and
the gauge invariant uniqueness result used to establish
 our results. Section \ref{os} deals with out-split graphs,
section \ref{dels} with in- and out-delays, and section \ref{is}
with in-split graphs. Finally section \ref{shifteq} relates our
results to those in \cite{b1}. Our main results are:

\begin{itemize}
\item[1.] If $E$ is a directed graph and $F$ is a proper out-split
graph formed from $E$ then  $C^* (E) \cong C^* (F)$ (Theorem
\ref{outsplitiso}).

\item[2.] If $E$ is a directed graph and $F$ is an out-delayed graph formed
from $E$ then $C^* (E)$ is strongly Morita equivalent to $C^* (F)$
if and only if $F$ arises from a proper out-delay (Theorem
\ref{outdelaysme}).

\item [3.] If $E$ is a directed graph and $F$ is an in-delayed
graph formed from $E$ then $C^* (E)$ is strongly Morita equivalent
to $C^* (F)$ (Theorem \ref{indelayme}).

\item[4.] If $E$ is a directed graph and $F$ is an in-split graph
formed from $E$ then $C^* (E)$ is strongly Morita equivalent to
$C^* (F)$ if and only if $F$ arises from a proper in-splitting
(Corollary \ref{insplitme}).
\end{itemize}

\section{The $C^*$-algebra of a directed graph} \label{basic}

Here we briefly set out some of the basic definitions and
terminology which we use throughout this paper. A directed graph
$E$  consists of countable sets of vertices and edges $E^0$ and
$E^1$, together with maps $r, s : E^1 \to E^0$ giving direction of
each edge. The maps $r, s$ extend naturally to $E^*$, the
collection of all finite paths in $E$. The graph $E$ is called
{\em row-finite} if every vertex emits a finite number of edges.

A Cuntz-Krieger $E$-family consists of a collection $\{ s_e : e
\in E^1 \}$ of partial isometries with orthogonal ranges, and
mutually orthogonal projections $\{ p_v : v \in E^0 \}$ satisfying
\begin{itemize}
\item[(i)] $s_e^* s_e = p_{r(e)}$,

\item[(ii)] $s_e s_e^* \le p_{s(e)}$,

\item[(iii)] if $v$ emits finitely many edges then $p_v =
\sum_{s(e)=v} s_e s_e^*$.

\end{itemize}

\noindent The graph $C^*$-algebra of $E$, $C^* (E)$ is the
universal $C^*$-algebra generated by a Cuntz-Krieger $E$-family.
An important property of a graph $C^*$-algebra is the
existence of an action $\gamma$ of ${\bf T}$, called the gauge
action, which is characterised by
\[
\gamma_z s_e = z s_e , \text{ and } \gamma_z p_v = p_v
\]

\noindent where $\{s_e , p_v \}\subseteq C^* (E)$ is the canonical
Cuntz-Krieger $E$-family and $z \in {\bf T}$. This gauge action is
a key ingredient in the uniqueness theorem which we shall
frequently use:

\begin{thm} \label{bhrsz} \cite[Theorem 2.1]{bhrsz}
Let $E$ be a directed graph, $\{ S_e , P_v \}$ be a Cuntz-Krieger
$E$-family and $\pi : C^*(E) \to C^*( S_e,P_v )$ the homomorphism
satisfying $\pi (s_e) = S_e$ and $\pi (p_v) = P_v$. Suppose that
each $P_v$ is non-zero, and that there is a strongly continuous
action $\beta$ of ${\bf T}$ on $C^* ( S_e , P_v )$ such that
$\beta_z \circ \pi = \pi \circ \gamma_z$ for all $z \in {\bf T}$.
Then $\pi$ is faithful.
\end{thm}

\noindent To apply Theorem \ref{bhrsz}, we exhibit a non-trivial
Cuntz-Krieger $E$-family within a $C^*$-algebra $B$, which carries
a suitable ${\bf T}$-action $\beta$.

Some results in this paper require the following result on Morita
equivalence of graph algebras. As in \cite[Remark 3.1]{bhrsz}  we
define the saturation $\Sigma H(S)$ of $S \subseteq E^0$ to be the
union of the sequence $\Sigma_n(S)$ of subsets of $E^0$ defined
inductively as follows:
\begin{eqnarray*}
\Sigma_0(S) &:= & \{ v \in E^0 \;:\; v = r(\mu) \mbox{ for some }
\mu \in E^* \mbox{ with } s(\mu) \in S \}
\\
\Sigma_{n+1}(S)&:=& \Sigma_n(S) \cup \{w \in E^0 \; : \; 0 <
|s^{-1}(w)| < \infty \mbox{ and } s(e) = w \mbox{ imply } r(e) \in
\Sigma_n(S)\}
\end{eqnarray*}
We note that if $E$ is row-finite then $\Sigma H(S)$ is the
saturation of the hereditary set $\Sigma_0(S)$ as defined in
\cite{bprs}.

\begin{lem} \label{1.3}
Suppose that $E$ is a directed graph, $S$ a subset of $E^0$ and
$\{s_e,p_v\}$ the canonical Cuntz-Krieger $E$-family. Let $P =
\sum_{v \in S}p_v$.  Then $P \in {\mathcal M}(C^*(E))$ and the
corner $P C^*(E) P$ is full if and only if $\Sigma H (S) = E^0$.
\end{lem}

\begin{proof}
By \cite[Lemma 3.3.1]{pr} the sum $\sum_{v \in S}p_v$ converges to
a projection $P \in {\mathcal M}(C^*(E))$. We claim that $PC^*(E)P
\subseteq I_{\Sigma H(S)}$.  Let $s_{\mu} s_{\nu}^*$ be a nonzero
element of $PC^*(E)P$, then $s(\mu) \in S$ and so $p_{s(\mu)} \in
I_{\Sigma H (S)}$.  Thus $s_{\mu} s_{\nu}^* = p_{s(\mu)} s_{\mu}
s_{\nu}^* \in I_{\Sigma H (S)}$
establishing our claim.

If $PC^*(E)P$ is full then $I_{\Sigma H(S)} = C^* (E)$ and so
$\Sigma H(S) = E^0$ by \cite[Section 3]{bhrsz}. Conversely,
suppose that $\Sigma H(S) = E^0$ and $PC^*(E)P \subseteq I$ for
some ideal $I$ in $C^*(E)$. By \cite[Lemma 3.2]{bhrsz} $H_I = \{v
: p_v \in I \}$ is a saturated hereditary subset of $E^0$
containing $S$ and hence $\Sigma H(S)$. Thus $C^* (E) = I_{\Sigma
H(S)} \subseteq I$ and the result follows.
\end{proof}

\section{Out-Splittings} \label{os}

The following definitions are adapted from \cite[Definition
2.4.3]{lm}. Let $E = ( E^0 , E^1 , r , s )$ be a directed graph.
For each $v \in E^0$ which emits an edge, partition $s^{-1} (v)$
into disjoint nonempty subsets $\mathcal{E}^1_v , \ldots ,
\mathcal{E}^{m(v)}_v$ where $m(v) \ge 1$ (if $v$ is a sink, then
we put $m(v)=0$). Let $\mathcal{P}$ denote the resulting partition
of $E^1$. We form the {\em out-split graph $E_s ( \mathcal{P} )$
from $E$ using $\mathcal{P}$} as follows: Let
\begin{align*}
E_s ( \mathcal {P} )^0 &= \{ v^i : v \in E^0 , 1 \le i \le m(v) \} \cup \{ v :
m(v)=0 \} , \\
E_s ( \mathcal{P} )^1 &= \{ e^j : e \in E^1, 1 \le j \le m ( r (e) )
\} \cup \{ e : m ( r(e) ) = 0 \} ,
\end{align*}

\noindent
and define $r_{E_s ( \mathcal{P} )} , s_{E_s ( \mathcal{P} )} : E_s (
\mathcal{P} )^1 \rightarrow E_s ( \mathcal{P} )^0$ for $e \in
\mathcal{E}^i_{s(e)}$ by
\begin{align*}
s_{E_s ( \mathcal{P} )} ( e^j ) &= s(e)^i \text{ and } s_{E_s (
\mathcal{P} )} (e)
= s(e)^i  \\
r_{E_s ( \mathcal{P} )} ( e^j ) &= r(e)^j \text{ and }  r_{E_s (
\mathcal{P} )} ( e ) = r(e) .
\end{align*}

\noindent The partition $\mathcal{P}$ is {\em proper} if for every
vertex $v$ with infinite valency we have $m(v) < \infty$ and only
one of the partition sets $\mathcal{E}_v^i$ is infinite.

\begin{examples} \label{firstex}
\begin{itemize}
\item[(i)] The partitions which give rise to the out-splittings
described in \cite[Figures 2.4.3--2.4.5]{lm} and \cite[\S 4.1]{d}
are all examples of proper partitions.

\item[(ii)] If we out-split at an infinite valence vertex, taking
a partition $\mathcal{P}$ which has finitely many subsets, such as
in
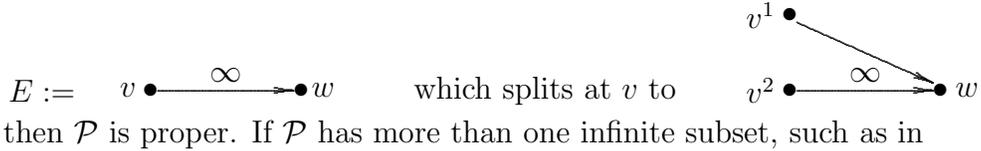
\[
\beginpicture
\setcoordinatesystem units <1cm,1cm>

\setplotarea x from -0.5 to 11, y from 0 to 1

\put{$E:=$}[r] at -0.5 0

\put{$\bullet$} at 0.5 0

\put{$\bullet$} at 2.5 0

\put{$v$} at 0.2 0

\put{$w$} at 2.8 0

\put{$\infty$} at 1.5 0.2

\arrow <0.25cm> [0.1,0.3] from 0.6 0 to  2.4 0

\put{which splits at $v$ to}[l] at 4 0

\put{$\bullet$} at 9 1

\put{$\bullet$} at 9 0

\put{$\bullet$} at 11 0

\put{$v^1$}[r] at 8.8 1

\put{$v^2$}[r] at 8.8 0

\put{$w$}[l] at 11.2 0

\put{$\infty$} at 10 0.2

\arrow <0.25cm> [0.1,0.3] from 9.1 0  to  10.9 0

\arrow <0.25cm> [0.1,0.3] from 9.1 0.9 to  10.9 0.1

\endpicture
\]

\noindent then $\mathcal{P}$ is proper. If $\mathcal{P}$ has more
than one infinite subset, such as in
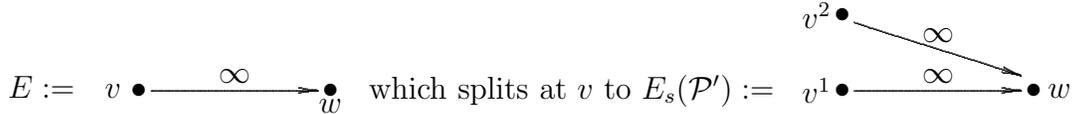
\[
\beginpicture
\setcoordinatesystem units <1.7cm,1cm>

\setplotarea x from 0 to 7, y from -0.2 to 1

\put{$E:=$}[r] at -0.5 0.05

\put{$\bullet$} at 0 0

\put{$\bullet$} at 1.5 0

\put{$\infty$} at 0.75 0.2

\put{$v$} at -0.2 0

\put{$w$} at 1.5 -0.2

\arrow <0.25cm> [0.1,0.3] from 0.1 0 to  1.4 0

\put{which splits at $v$ to $E_s ( \mathcal{P}'):=$}[l] at 1.8 0

\put{$\bullet$} at 5.5 0

\put{$\bullet$} at 5.5 1

\put{$\bullet$} at 7 0

\put{$v^1$} at 5.3 0

\put{$v^2$} at 5.3 1

\put{$w$} at 7.2 0

\put{$\infty$} at 6.25 0.2

\put{$\infty$}[b] at 6.25 0.65

\arrow <0.25cm> [0.1,0.3] from 5.6 0 to  6.9 0

\arrow <0.25cm> [0.1,0.3] from 5.6 0.9 to  6.9 0.2

\endpicture
\]

\noindent then $\mathcal{P}$ is not proper. Again $C^*(E)$ is not
Morita equivalent to $C^*(E_s({\mathcal P}))$ since the latter has
an additional ideal. If $\mathcal{P}$ has infinitely many subsets,
such as in
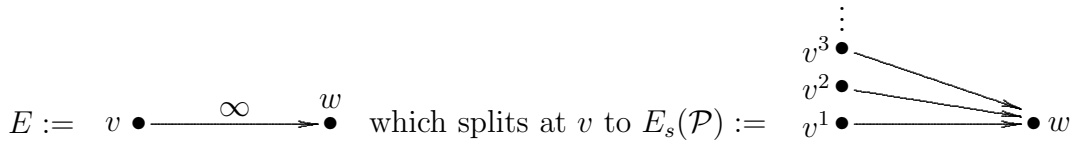
\[
\beginpicture
\setcoordinatesystem units <1.7cm,0.5cm>

\setplotarea x from 0 to 7, y from 0 to 3

\put{$E:=$}[r] at -0.5 0.05

\put{$\bullet$} at 0 0

\put{$\bullet$} at 1.5 0

\put{$\infty$} at 0.75 0.4

\put{$v$} at -0.2 0

\put{$w$} at 1.5 0.6

\arrow <0.25cm> [0.1,0.3] from 0.1 0 to  1.4 0

\put{which splits at $v$ to $E_s ( \mathcal{P}):=$}[l] at 1.8 0

\put{$\bullet$} at 5.5 0

\put{$\bullet$} at 5.5 1

\put{$\bullet$} at 5.5 2

\put{$\vdots$} at 5.5 3

\put{$\bullet$} at 7 0

\put{$v^1$} at 5.3 0

\put{$v^2$} at 5.3 1

\put{$v^3$} at 5.3 2

\put{$w$} at 7.2 0

\arrow <0.25cm> [0.1,0.3] from 5.6 0 to  6.9 0

\arrow <0.25cm> [0.1,0.3] from 5.6 0.9 to  6.9 0.2

\arrow <0.25cm> [0.1,0.3] from 5.6 1.9 to  6.9 0.4

\endpicture
\]

\noindent then $\mathcal{P}$ is not proper.  In this case,
$C^*(E)$ and $C^*(E_s({\mathcal P}))$ are not Morita equivalent:
The former is non-simple and the latter simple.
\end{itemize}
\end{examples}

\begin{thm} \label{outsplitiso}
Let $E$ be a directed graph, $\mathcal{P}$ a partition of $E^1$
and $E_s ( \mathcal{P} )$ the out-split graph formed from $E$
using $\mathcal{P}$. If $\mathcal{P}$ is proper then $C^* (E)
\cong C^* ( E_s ( \mathcal{P} ) )$.
\end{thm}

\begin{proof}
Let $\{ s_f , p_w : f \in E_s ( \mathcal{P} )^1 , w \in E_s (
\mathcal{P} )^0 \}$ be a Cuntz-Krieger $E_s ( \mathcal{P}
)$-family. For $v \in E^0$ and $e \in E^1$
set $Q_v = p_v$ if $m(v)=0$, $T_e = s_e$ if $m(r(e))=0$,
\[
Q_v = \sum_{1 \le i \le m(v)} p_{v^i} \text{ if } m(v) \neq 0 \text{ and } T_e =
\sum_{1 \le j \le  m ( r (e ) )} s_{e^j}  \text{ if } m(r(e)) \neq 0 .
\]

\noindent Because $\mathcal{P}$ is proper $m (v) < \infty$ for all
$v \in E^0$ and all of these sums are finite. We claim that $\{
T_e , Q_v : e \in E^1 , v \in E^0 \}$ is a Cuntz-Krieger
$E$-family in $C^* ( E_s ( \mathcal{P} ) )$.

The $Q_v$'s are non-zero mutually orthogonal projections since
they are sums of projections satisfying the same properties. The
partial isometries $T_e$ for $e \in E^1$ have mutually orthogonal
ranges since they consist of sums of partial isometries with
mutually orthogonal ranges. For $e \in E^1$ it is easy to see that
$T_e^* T_e = Q_{r(e)}$ and $T_e T_e^* \le Q_{s(e)}$.

For $e \in E^1$ with $m(r(e)) \neq 0$, then since $r_{E_s} ( e^j )
\neq r_{E_s} ( e^k )$ for $j \neq k$ we have
\begin{equation} \label{reltwo}
T_e T_e^* = \Big( \sum_{1 \le j \le m ( r(e) )} s_{e^j} \Big)
\Big( \sum_{1 \le k \le m ( r(e) )} s_{e^k} \Big)^* = \sum_{1 \le
j \le m(r(e))} s_{e^j} s_{e^j}^* .
\end{equation}

\noindent If $m(r(e))=0$ then $T_e T_e^* = s_e s_e^*$. For $v \in
E^0$ and $1 \le i \le m(v)$ put $\mathcal{E}_{1,v}^i = \{ e \in
\mathcal{E}_v^i : m(r(e)) \geq 1 \}$ and $\mathcal{E}_{0,v}^i = \{
e \in \mathcal{E}_v^i : m(r(e)) = 0 \}$. If $v \in E^0$ has finite
valency and is not a sink then $s^{-1}(v) =
\bigcup_{i=1}^{m(v)}{\mathcal E}^i_v$ and for $1 \le i \le m(v)$
we have
\[
s_{E_s ( \mathcal{P} )}^{-1}
( v^i ) = \{ e^j : e \in \mathcal{E}_{1,v}^i , 1 \le j \le m(r(e)) \} \cup 
{\mathcal E}^i_{0,v}
\]

\noindent
Hence using (\ref{reltwo}) we may compute
\begin{align*}
Q_v = \sum_{1 \le i \le m(v)} p_{v^i} &= \sum_{1 \le i \le m(v)} \sum_{e \in
\mathcal{E}_{1,v}^i} \sum_{1 \le j \le m(r(e))} s_{e^j} s_{e^j}^* +
\sum_{1 \le i \le m(v)} \sum_{e \in
\mathcal{E}_{0,v}^i} s_e s_e^*  \\
&= \sum_{1 \le i \le m(v)} \; \sum_{e \in \mathcal{E}_v^i} T_e
T_e^* = \sum_{e: s(e)=v} T_e T_e^* ,
\end{align*}

\noindent completing the proof of our claim, since vertices $v \in
E^0$ with $m(v)=0$ are sinks.

Let $\{ t_e , q_v \}$ be the canonical generators of $C^* (E)$,
then by the universal property of $C^* (E)$ there is a
homomorphism $\pi : C^* (E) \rightarrow C^* ( E_s ( \mathcal{P}
))$ taking $t_e$ to $T_e$ and $q_v$ to $Q_v$. To prove that $\pi$
is onto we show that the generators of $C^* ( E_s ( \mathcal{P} )
)$ lie in $C^* ( T_e , Q_v )$. If $w = v^j \in E_s(\mathcal{P})^0$
has finite valency or $w$ is a sink, then $p_w \in C^*( T_e , Q_v
)$. If $v^j $ has infinite valency, then without loss of
generality we suppose $j=1$. Since $\mathcal{P}$ is proper it
follows that $v^2 , \ldots , v^{m(v)}$ have finite valency, so
$p_{v^2} , \ldots p_{v^{m(v)}} \in C^* ( T_e , Q_v )$ and hence
\[
p_{v^1} = Q_v - p_{v^2} - \ldots - p_{v^{m(v)}} \in C^* ( T_e ,
Q_v  ) .
\]

\noindent If $e^j \in E_s ( \mathcal{P} )^1$ then $m(r(e)) \neq
0$. Since $p_{r(e)^j} \in C^* ( T_e , Q_v )$ we have $s_{e^j} =
T_e p_{r(e)^j}
 \in C^* ( T_e , Q_v )$. If $e \in E_s ( \mathcal{P}
)^1$ then $m(r(e)) = 0$ and so $s_e = T_e \in C^* ( T_e , Q_v )$.

Since $\pi$ commutes with the canonical gauge action on each
$C^*$-algebra and as $Q_v \neq 0$ for all $v \in E^0$ it follows
from Theorem \ref{bhrsz} that $\pi$ is injective, and the result
follows.
\end{proof}

\begin{rmks} \label{osrem}
\begin{itemize}
\item[(i)] The {\em maximal out-splitting} $\tilde{E}$ of $E$ is
formed from a partition $\mathcal{P}$ of $E^1$ which admits no
refinements. For a graph $E = ( E^0 , E^1 , r , s)$ without sinks
$\tilde{E}$ is isomorphic to the dual graph $\widehat{E} = ( E^1 ,
E^2 , r' , s' )$  (where $r'(ef)=f$ and $s'(ef)=e$): Since the
out-splitting is maximal and there are no sinks, we have
\[
\tilde{E}^0 = \{ v^e : s(e) = v \} , \text{ and }
\tilde{E}^1 = \{ e^f : s(f)=r(e) \} .
\]

\noindent The maps $v^e \mapsto e$ and $e^f \mapsto ef$ induce an
isomorphism from $\tilde{E}$ to $\widehat{E}$. We define the dual
graph $\widehat{E}$ of any directed graph $E$ to be its maximal
out-split graph $\tilde{E}$. Since a maximal out-splitting is
proper if and only $E$ is row-finite we may now use Theorem
\ref{outsplitiso} to generalise \cite[Corollary 2.5]{bprs} to any
row-finite graph.

\item[(ii)] Brenken defines a graph algebra $G^*(E)$ which, under
certain conditions, is isomorphic to $C^*(E)$.  By \cite[Theorem
3.4]{bb}, $G^* (E)$ and the $C^*$-algebra of its out-splitting are
isomorphic if the graph satisfies a certain properness condition
(see \cite[Definition 3.2]{bb}). This result only reduces to $C^*
(E)$ when $E$ is row-finite.

\item[(iii)] For row-finite graphs with finitely many vertices and
no sinks, a proof of Theorem \ref{outsplitiso} may be deduced from
\cite[\S 4.1]{d}. In \cite[Section 6]{ds} similar results are
proved for an ``explosion" which is a particular example of an
out-splitting operation.

\item[(iv)] The $C^*$-algebra of an out-split graph is isomorphic
to the $C^*$-algebra of an ultragraph (see \cite{t}). Given a
directed graph $E$ and a partition $\mathcal{P}$ define the
ultragraph $\mathcal{G} ( \mathcal{P} ) = ( G^0 , \mathcal{G}^1 ,
r' , s' )$ as follows: Put $G^0 = E_s ( \mathcal{P} )^0$,
$\mathcal{G}^1 = E^1$, $s' ( e ) = s(e)^i$ if $e \in
\mathcal{E}_{s(e)}^i$, $r' (e) = r(e)$ if $m(r(e))=0$ and
\[
r' ( e ) = \{ r(e)^i : 1 \le i \le m(r(e)) \} \text{ if } m(r(e)) \neq 0 .
\]

\noindent We claim that if $\mathcal{P}$ is proper then $C^* (
\mathcal{G} ( \mathcal{P} ) ) \cong C^* ( E_s ( \mathcal{P} ) )$.
When $\mathcal{P}$ is proper $\mathcal{G}^0$ is the set of finite
subsets of $E_s ( \mathcal{P} )^0$. Let $\{ p_A , s_e : A \in
\mathcal{G}^0 , e \in \mathcal{G}^1 \}$ be a Cuntz-Krieger
$\mathcal{G} ( \mathcal{P} )$-family. For $w \in E_s ( \mathcal{P}
)^0$ set $Q_{w} = p_{w}$ and for $f^j \in E_s ( \mathcal{P} )^1$
put $T_{f^j} = s_f p_{r(f)^j}$.  Then $\{ Q_w , T_{f^j} \}$ is a
Cuntz-Krieger $E_s ( \mathcal{P} )$-family in which each $Q_w \ne
0$. Let $\{ t_{f^j} , q_w \}$ be the canonical generators of $C^*
( E_s ( \mathcal{P} ) )$, then by the universal property of $C^*
(E_s ( \mathcal{P} ) )$ there is a map $\pi : C^* ( E_s (
\mathcal{P} ) ) \to C^* ( \mathcal{G} ( \mathcal{P} ) )$ sending
$t_{f^j}$ to $T_{f^j}$ and $q_w$ to $Q_w$. As each $A \in
\mathcal{G}^0$ is finite, $C^* ( T_{f^j} , Q_w )$ contains each
generator of $C^* ( \mathcal{G} ( \mathcal{P} ) )$, so $\pi$ is
onto. By \cite[Section 2]{t} there is an appropriate action of
${\bf T}$ on $C^*(\mathcal{G} (\mathcal{P}))$, so $\pi$ is
injective by Theorem \ref{bhrsz}, proving our claim.

\item[(v)] Let $\Gamma$ act freely on the edges of a row-finite
graph $E$, then the induced $\Gamma$-action on the dual graph
$\widehat{E}$ is free on its vertices. The isomorphism $C^* (E)
\cong C^* ( \widehat{E} )$ is $\Gamma$-equivariant, so $C^* ( E )
\times \Gamma \cong C^* ( \widehat{E} ) \times \Gamma$ and by
\cite[Corollary 3.3]{kqr} we have
\begin{equation} \label{obs}
C^* (E) \times \Gamma \cong C^* ( \widehat{E} / \Gamma ) \otimes
\mathcal{K} \left( \ell^2 ( {\bf Z}_2 ) \right).
\end{equation}

\noindent For instance, there is a free action of ${\bf Z}_2$ on
the edges of graph $B_2$ which consists of a single vertex and two
edges (the ``flip-flop automorphism'' of $C^* ( B_2 ) \cong
\mathcal{O}_2$ described by \cite{ar}). Equation (\ref{obs}) shows
that $\mathcal{O}_2 \times {\bf Z}_2$ is Morita equivalent
$\mathcal{O}_2$.
\end{itemize}
\end{rmks}

\section{Delays} \label{dels}

Let $E = ( E^0, E^1, r, s )$ be a directed graph. A map $d_s : E^0
\cup E^1 \rightarrow {\bf N} \cup \{ \infty \}$ such that
\begin{itemize}
\item[(i)] if $w \in E^0$ is not a sink then $d_s (w) = \sup \{ d_s (e) : s(e)=w \}$,
\item[(ii)] if $d_s ( x ) = \infty$ for some $x$ then either $x$ is a sink or $x$
emits infinitely many edges
\end{itemize}

\noindent
is called a {\em Drinen source-vector}. Note that only vertices are
allowed to have an infinite $d_s$-value; moreover if $d_s(v) = \infty$
and $v$ is not a sink, then there
are edges with source $v$ and arbitrarily large $d_s$-value. From this
data we construct a new graph as follows: Let
\begin{align*}
d_s (E)^0 &= \{ v^i : v \in E^0 , 0 \le i \le d_s (v) \} \text{ and } \\
d_s (E)^1 &= E^1 \cup \{ f(v)^i : 1 \le i \le d_s (v)\} ,
\end{align*}

\noindent
and for $e \in E^1$ define $r_{od} (e) = r(e)^0$ and $s_{od} (e) = s(e)^{d_s
(e)}$. For $f (v)^i$ define $s_{od} ( f(v)^{i} ) = v^{i-1}$ and
$r_{od} ( f(v)^i ) = v^i$.  The resulting directed graph $d_s (E)$ is
called the {\em out-delayed graph of $E$} for the Drinen source-vector
$d_s$.

In the out-delayed graph the original vertices correspond to those
vertices with superscript $0$; the edge $e \in E^1$ is delayed
from leaving $s(e)^0$ and arriving at $r(e)^0$ by a path of length
$d_s (e)$. The Drinen source vector $d_s$ is {\em strictly proper}
if, whenever $v$ has infinite valency, there is no $v^i$ with
infinite valency unless $i=d_s (v) < \infty$. A Drinen source
vector $d_s$ which gives rise to an out-delayed graph $d_s (E)$
which may be constructed using a finite sequence of strictly
proper Drinen source vectors is said to be {\em proper}.

\begin{examples} \label{outdelayex}
\begin{itemize}

\item[(i)] The notion of an out-delay in the context of graph
$C^*$-algebras was first introduced in \cite[\S 4]{ck} and
subsequently generalised in \cite{d}. The graphs shown in \cite[\S
3.1]{d} are all examples of out-delays for some proper Drinen
source-vector where all the edges out of a given vertex are
delayed by the same amount.

\item[(ii)] The Drinen-Tomforde desingularisation of a graph
described in \cite[Definition 2.2]{dt} is an example of a
out-delay with a proper Drinen source-vector: If $v$ has infinite
valency then the edges with source $v$ may be written as $\{ e_i :
i \in {\bf N} \}$; we set $d_s (v) = \infty$ and $d_s (e_i ) = i$
for $i \in {\bf N}$. If $v$ has finite valency then we set $d_s
(v)=0$ (and so $d_s (e)=0$ for all $e \in s^{-1} (v)$). If $v$ is
a sink then we put $d_s (v) = \infty$. The resulting graph $d_s
(E)$ is row-finite with no sinks.

\item[(iii)] Putting $d_s (v) = \infty$ for a sink adds an
infinite tail to the sink. If $d_s (v) = \infty$ for all sinks and
$d_s (v) = 0$ for all vertices which emit edges then $d_s (E)$ is
the graph $E$ with tails added to all sinks (cf.\ \cite[Lemma
1.4]{rs}).

\item[(iv)] Consider the graph $E$ shown below.
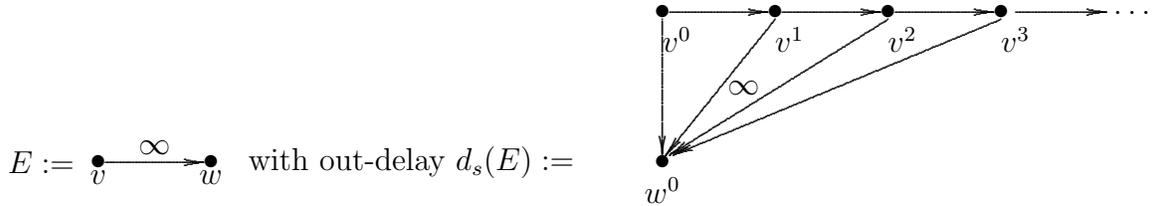
\[
\beginpicture
\setcoordinatesystem units <1cm,1cm>

\setplotarea x from -1.2 to 14, y from -0.2 to 2.2

\put{$E:=$}[l] at -1.2 0

\put{$\bullet$} at 0 0

\put{$\bullet$} at 1.5 0

\put{$\infty$} at 0.75 0.2

\put{$v$} at 0 -0.2

\put{$w$} at 1.5 -0.2

\arrow <0.25cm> [0.1,0.3] from 0.1 0 to 1.4 0

\put{with out-delay $d_s (E):=$}[l] at 2 0

\put{$\ldots$}[l] at 13.5 2

\put{$\bullet$} at 7.5 2

\put{$\bullet$} at 9 2

\put{$\bullet$} at 10.5 2

\put{$\bullet$} at 12 2

\put{$\bullet$} at 7.5 0

\put{$v^3$}[tl] at 12 1.8

\put{$v^2$}[tl] at 10.5 1.8

\put{$v^1$}[tl] at 9 1.8

\put{$v^0$}[tl] at 7.5 1.8

\put{$w^0$}[t] at 7.5 -0.2

\put{$\infty$}[r] at 8.79 1

\arrow <0.25cm> [0.1,0.3] from 7.6 2 to 8.9 2

\arrow <0.25cm> [0.1,0.3] from 9.1 2 to 10.4 2

\arrow <0.25cm> [0.1,0.3] from 10.6 2 to 11.9 2

\arrow <0.25cm> [0.1,0.3] from 12.2 2 to 13.4 2

\arrow <0.25cm> [0.1,0.3] from 12 1.9 to 7.65 0.1

\arrow <0.25cm> [0.1,0.3] from 10.5 1.9 to 7.6 0.1

\arrow <0.25cm> [0.1,0.3] from 9 1.9 to 7.55 0.1

\arrow <0.25cm> [0.1,0.3] from 7.5 1.9 to 7.5 0.1

\endpicture
\]

\noindent The Drinen source-vector for this out-delay is not
proper since vertex $v_1$ has infinite valency. Moreover,  the
$C^*$-algebra $C^* ( d_s(E))$ is not Morita equivalent to $C^*
(E)$ since the former $C^*$-algebra has $2$ proper ideals and the
latter only one.
\end{itemize}
\end{examples}

\begin{thm}\label{outdelaysme}
Let $E$ be a directed graph and $d_s : E^0 \cup E^1 \rightarrow
{\bf N} \cup \{ \infty \}$ be a Drinen source-vector. Then $C^* (
d_s ( E ) )$ is strongly Morita equivalent to $C^* (E)$ if and
only if $d_s$ is proper.
\end{thm}

\begin{proof}
Without loss of generality we may assume that $d_s : E^0 \cup E^1
\rightarrow {\bf N} \cup \{ \infty \}$ is an essentially proper
Drinen source-vector. Let $\{ s_f , p_w : f \in d_s (E)^1 , w \in
d_s (E)^0 \}$ be a Cuntz-Krieger $d_s (E)$-family. For $e \in E^1$
and $v \in E^0$ define $Q_v = p_{v^0}$ and
\[
T_e = s_{f(s(e))^1} \ldots s_{f(s(e))^{d_s (e)}} s_e \text{ if } d_s (e) \neq 0
\text{ and } T_e = s_e \text{ otherwise.}
\]

\noindent We claim that $\{ T_e , Q_v : e \in E^1 , v \in E^0 \}$
is a Cuntz-Krieger $E$-family. The $Q_v$'s are non-zero mutually
orthogonal projections since the $p_v$'s are. The $T_e$'s are
partial isometries with mutually orthogonal ranges since they are
products of partial isometries with this property. For $e \in E^1$
it is routine to check that $T_e^* T_e = Q_{r(e)}$ and $T_e T_e^*
\le Q_{s(e)}$.

If $v \in E^0$ is neither a sink nor has infinite valence, then
$d_s (v) < \infty$. If $d_s (v) =0$, then we certainly have $Q_v =
\sum_{s(e)=v} T_e T_e^*$. Otherwise, for $0 \le j \le d_s (v)-1$
we have
\begin{equation} \label{gantletrels}
p_{v^j} = \sum_{s(e)=v , d_s (e)=j} s_e s_e^* + s_{f(v)^{j+1}} p_{v^{j+1}}
s_{f(v)^{j+1}}^* ,
\end{equation}

\noindent
and since we must have some edges with $s(e)=v$ and $d_s (e) = d_s (v)$ we have
\begin{equation} \label{endcond}
p_{v^{d_s (v)}} = \sum_{s(e)=v , d_s (e) = d_s (v)} s_e s_e^* .
\end{equation}

\noindent Using (\ref{gantletrels}) recursively and
(\ref{endcond}) when $j=d_s (v)-1$ we see that
\[
Q_v = p_{v^0} = \sum_{s(e)=v,d_s(e)=0} T_e T_e^* + \ldots +
\sum_{s(e)=v,d_s(e)=d_s(v)} T_e T_e^* = \sum_{s(e)=v} T_e T_e^* ,
\]
and this establishes our claim.

Let $\{ t_e , q_v \}$ be the canonical generators of $C^* (E)$,
then by the universal property of $C^* (E)$ there is a
homomorphism $\pi : C^* (E) \rightarrow C^* ( d_s ( E ) )$ which
takes $t_e$ to $T_e$ and $q_v$ to $Q_v$. It remains to show that
$C^*( T_e, Q_v )$ is a full corner in $C^*(d_s(E))$.

Let $\alpha$ denote the strongly continuous ${\bf T}$-action
satisfying, for $z \in {\bf T}$,
\[
\alpha_z s_e = z s_e , ~ \alpha_z s_{f(v)^i} = s_{f(v)^i} \text{
for } 1 \le i \le d_s(v) \text{ and } \alpha_z p_{v^i} = p_{v^i}
\text{ for } 0 \le i \le d_s(v).
\]

\noindent It is straightforward to check that $\pi \circ \gamma =
\alpha \circ \pi$ where $\gamma$ is the usual gauge action of
${\bf T}$ on $C^* (E)$ and it follows from Theorem \ref{bhrsz}
that $\pi$ is injective.

By Lemma \ref{1.3} the sum $\sum_{v \in E^0}p_{v^0}$ converges to
a projection $P \in {\mathcal M}(C^*(d_s(E)))$.  We claim that
$C^*( T_e, Q_v )$ is equal to $P C^*( d_s(E)) P$. Note that if $P
s_{\mu} s_{\nu}^* P = s_{\mu} s_{\nu}^* \neq 0 $ then $s_{od} (
\mu ) = v^0$ and $s_{od} ( \nu ) = w^0$ for some $v, w \in E^0$,
and $r_{od} ( \mu ) = r_{od} (\nu)$.

If $\mu = \nu = v^0$ then $s_\mu s_\nu^* = p_{v^0} = Q_v \in C^*(
\{T_e, Q_v \} )$. If $r_{od} (\mu) = u^0$ for some $u \in E^0$
then there are paths $\alpha, \beta \in E^*$ with $r(\alpha) =
r(\beta) = u$ such that $s_{\mu}s_{\nu}^* = T_{\alpha}
T_{\beta}^*$  and so $s_{\mu}s_{\nu}^* \in C^*( T_e,Q_v )$.
Suppose now that $r_{od} ( \mu ) \notin E^0$. Then $r_{od}( \mu )
= u^q$ for some $u \in E^0$ and $0 < q \leq d_s (u)$ and we can
write
\[
s_\mu s_\nu^* = T_\alpha s_{f(u)^1} \ldots s_{f(u)^{q}} s_{f(u)^{q}}^* \ldots s_{f(u)^1}^* T_\beta^* .
\]
\noindent for some $\alpha, \beta \in E^*$. Suppose $q > 1$, then
since $d_s$ is proper, $f(u)^{q-1}$ has finite valency. If there
are no edges in $E$ with $s(e)=u$ and $d_s (e)=q-1$ then
\begin{equation} \label{one}
s_\mu s_\nu^* = T_\alpha s_{f(u)^1} \ldots s_{f(u)^{q-1}} s_{f(u)^{q-1}}^* \ldots s_{f(u)^1}^* T_\beta^* .
\end{equation}

\noindent If there are a finite number of edges $e_1 , \ldots ,
e_l \in E^1$ with $s ( e_i ) = u$ and $d_s ( e_i )=q-1$ for $i= 1
, \ldots , l$, then
\begin{equation} \label{two}
\begin{array}{rcl}
s_{\mu}s_{\nu}^* &=& T_{\alpha}s_{f(u)^1} \ldots s_{f(u)^{q-1}}(
p_{u^{q-1}} - s_{e_1} s_{e_1}^* - \ldots - s_{e_l} s_{e_l}^*)
s_{f(u)^{q-1}}^* \ldots s_{f(u)^1}^* T_{\beta}^* \\
&=& T_{\alpha} s_{f(u)^1} \ldots s_{f(u)^{q-1}} s_{f(u)^{q-1}}^*\ldots
s_{f(u)^1}^* T_{\beta}^* -
\sum_{i=1}^l T_{\alpha} T_{e_i} T_{e_i}^* T_{\beta}^*.
\end{array}
\end{equation}

\noindent Our new expression for $s_\mu s_\nu^*$ may now be
analysed as in (\ref{one}) or (\ref{two}), reducing the value of
$q$ until all $s_{f(u)^i}$ terms are removed. Then
\[
s_{\mu} s_{\nu}^* = T_\alpha T_\beta^* - \sum_{s(e)=u, d_s (e) \le
q-1} T_{\alpha e} T_{\beta e}^* \in C^* (  T_e , Q_v ) ,
\]

\noindent completing the proof of our claim.

Since $\Sigma H ( \{v^0: v \in E^0 \} ) = d_s(E)^0$, it follows
from Lemma \ref{1.3} that $P C^*(d_s(E)) P$ is a full corner in
$C^*(d_s(E))$ and hence $C^*(d_s(E))$ and $C^*(E)$ are strongly
Morita equivalent.

If $d_s$ is not proper, then there are at least two vertices
$f(v)^i , f(v)^j$ with $0 \le i \le j \le d_s (v)$ emitting
infinitely many edges. In this case there is an ideal generated by
$p_{f(v)^i}$ in $C^* ( d_s (E) )$ which was not present in $C^*
(E)$.
\end{proof}

\noindent We are grateful to Daniel Gow and Tyrone Crisp for
pointing out an error in an earlier version of Theorem
\ref{outdelaysme} (see also \cite{cg}).

\begin{rmks} \label{earlier}
\begin{itemize}
\item[(i)] Theorem \ref{outdelaysme} significantly generalises the
results in \cite[\S 3.1]{d}. Drinen shows a limited number of
graph groupoid isomorphisms for row-finite graphs with finitely
many vertices and no sinks in which each edge is equally delayed.

\item[(ii)] The desingularisation of a non row-finite graph is an
example of an out-delay (see Examples \ref{outdelayex} (ii)).
Moreover, any out-delay of a non row-finite graph using a proper
Drinen source vector with $d_s (v) = \infty$ for all vertices of
infinite valency provides an example of a row-finite graph $d_s
(E)$ whose $C^*$-algebra is Morita equivalent to $C^* (E)$. It
follows by \cite[Corollary 4.6]{ba2} that if $E$ satisfies
condition (K) then $\mathop{Prim} ( C^*  (E))$ is the primitive
ideal space of some AF-algebra.

\item[(iii)]  If $d_s : E^0 \cup E^1 \to {\bf N} \cup \{ \infty
\}$ is a Drinen source-vector then  $E$ is a deformation retract
of $d_s (E)$ (see \cite[\S 3.3]{st}). The construction of an
out-delayed graph replaces each vertex $v$ with $d_s (v) \ge 1$ by
the tree $\{ v^i : 0 \le i \le d_s (v) , f(v)^i : 1 \le i \le d_s
(v) \}$ which may be contracted to the root $v^0$ and identified
with $v$ (see also \cite[\S 1.5.5]{gt}). In particular $\pi_1(
d_s(E) ) \cong \pi_1(E)$ and the universal covering tree $T$ of
$E$ is a deformation retract of the universal covering tree $T'$
of $d_s (E)$. It follows that the boundary $\partial T$ of $T$ is
homeomorphic to the boundary $\partial T'$ of $T'$ (see \cite[\S
4]{kp}). Hence the Morita equivalence between $C^* (E)$ and $C^* (
d_s (E) )$ could be obtained for row-finite graphs with no sinks
using the Kumjian-Pask description of $C^*(E)$ as a crossed
product of $C_0 ( \partial T )$ by $\pi_1 (E)$ (see
\cite[Corollary 4.14]{kp}).
\end{itemize}
\end{rmks}

\noindent We now turn our attention to in-delays where edges are
delayed from arriving at their range. Let $E = ( E^0, E^1, r, s )$
be a graph. A map $d_r : E^0 \cup E^1 \rightarrow {\bf N} \cup \{
\infty \}$ satisfying
\begin{itemize}
\item[(i)] if $w$ is not a source then $d_r (w) = \sup \, \{ d_r
(e) : r(e)=w \}$,

\item[(ii)] if $d_r ( x ) = \infty$  then $x$ is either a source
or receives infinitely many edges
\end{itemize}

\noindent is called a {\em Drinen range-vector}. We construct a
new graph $d_r(E)$ called the {\em in-delayed graph of $E$} for
the Drinen range-vector $d_r$ as follows:
\begin{align*}
d_r (E)^0 &= \{ v_i : v \in E^0 , 0 \le i \le d_r (v) \} \text{ and } \\
d_r (E)^1 &= E^1 \cup \{ f(v)_i : 1 \le i \le d_r (v) \} ,
\end{align*}

\noindent and for $e \in E^1$ we define $r_{id} (e) = r(e)_{d_r
(e)}$ and $s_{id} (e) = s(e)_0$. For $f (v)_i$ we define $s_{id} (
f(v)_i ) = v_{i}$ and $r_{id} ( f(v)_i ) = v_{i-1}$.

\begin{examples} \label{indelayexample}
\begin{itemize}
\item[(i)] Consider the graph $E$ shown below, with edges $\{ e_i
: i \ge 0 \}$ from $v$ to $w$. If we set $d_r ( e_i ) = i$, $d_r (
v ) = 0$ and $d_r (w ) = \infty$ then
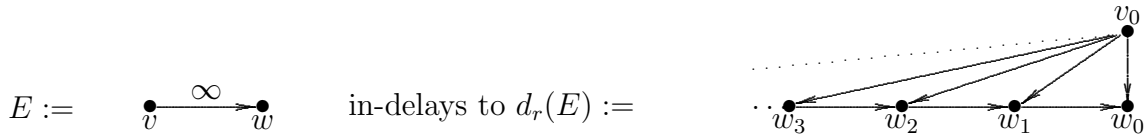
\[
\beginpicture
\setcoordinatesystem units <1cm,0.5cm>

\setplotarea x from -1 to 13, y from -0.2 to 2.2

\put{$E:=$}[r] at -1 0

\put{$\bullet$} at 0 0

\put{$\bullet$} at 1.5 0

\put{$\infty$} at 0.75 0.4

\put{$v$} at 0 -0.4

\put{$w$} at 1.5 -0.4

\arrow <0.25cm> [0.1,0.3] from 0.1 0 to 1.4 0

\put{ in-delays to $d_r (E):=$}[l] at 2.5 0

\put{$\ldots$}[l] at 8 0

\put{$\bullet$} at 8.5 0
\put{$\bullet$} at 10 0
\put{$\bullet$} at 11.5 0
\put{$\bullet$} at 13 0
\put{$\bullet$} at 13 2

\put{$w_0$}[t] at 13 -0.2
\put{$w_1$}[t] at 11.5 -0.2
\put{$w_2$}[t] at 10 -0.2
\put{$w_3$}[t] at 8.5 -0.2
\put{$v_0$}[b] at 13 2.2

\arrow <0.25cm> [0.1,0.3] from 8.6 0 to 9.9 0

\arrow <0.25cm> [0.1,0.3] from 10.1 0 to 11.4 0

\arrow <0.25cm> [0.1,0.3] from 11.6 0 to 12.9 0

\arrow <0.25cm> [0.1,0.3] from 13 1.9 to 13 0.1

\arrow <0.25cm> [0.1,0.3] from 12.9 1.9 to 11.6 0.1

\arrow <0.25cm> [0.1,0.3] from 12.85 1.9 to 10.1 0.1

\arrow <0.25cm> [0.1,0.3] from 12.8 1.9 to 8.6 0.1

\setdots

\setlinear

\plot 12.75 1.9 8 1 /

\endpicture
\]

\item[(ii)] Observe that putting $d_r (v) = \infty$ for a source
adds an infinite ``head'' to the source. If $d_r (v) = \infty$ for
all sources and $d_r (v)=0$ for all vertices which receive edges
then $d_r (E)$ is the graph $E$ with heads added to all sources
(cf. \cite[Lemma 1.4]{rs}).
\end{itemize}
\end{examples}

\begin{thm} \label{indelayme}
If $d_r : E^0 \rightarrow {\bf N} \cup \{
\infty \}$ is a Drinen range-vector, then $C^* ( d_r ( E ) )$ is strongly
Morita equivalent to $C^* (E)$.
\end{thm}

\begin{proof}
Let $d_r : E^0 \cup E^1 \rightarrow {\bf N} \cup \{ \infty \}$ be
a Drinen range-vector and $\{ s_e , p_v  : e \in d_r (E)^0 , v \in
d_r (E)^1 \}$ be a Cuntz-Krieger $d_r (E)$-family. For $v \in E^0$
let $Q_v = p_{v_0}$ and for $e \in E^1$ put
\[
T_e = s_e s_{f(r(e))_{d_r( e)}} \ldots s_{f(r(e))_{1}}  \text{ if } d_r (e)
\neq 0 \text{ and } T_e = s_e \text{ otherwise.}
\]

\noindent It is straightforward to check that $\{ T_e, Q_v \}$ is
a Cuntz-Krieger $E$-family in $C^*(d_r(E))$ in which all the
projections $Q_v$ are non-zero. Let $\{ t_e , q_v \}$ be the
canonical generators of $C^* (E)$, then by the universal property
of $C^* (E)$ there is a homomorphism $\pi : C^*(E) \to
C^*(d_r(E))$ satisfying $\pi ( t_e ) = T_e$ and $\pi  ( q_v ) =
Q_v$. It remains to show that $C^*( T_e, Q_v )$, the image $\pi$,
is a full corner in $C^*(d_r(E))$.

Let $\alpha$ be the strongly continuous ${\mathbf T}$-action
$\alpha$ on $C^*(d_r(E))$ satisfying, for $z \in {\bf T}$,
\[
\alpha_z(s_e) = zs_e , \alpha_z(s_{f(v)^i}) = s_{f(v)^i} \text{
for } 1 \leq i \leq d_r(v) \text{ and } \alpha_z(p_{v^i}) =
p_{v^i} \text{ for } 0 \leq i \leq d_r(v) . \]

\noindent It is straightforward to check that $\pi \circ \gamma =
\alpha \circ \pi$ where $\gamma$ is the usual gauge action on
$C^*(E)$ and it follows from Theorem \ref{bhrsz} that $\pi$ is
injective.

By Lemma \ref{1.3}, the sum $\sum_{v \in E^0} p_{v_0}$ converges
to a projection $P \in {\mathcal M}(C^*(d_r(E)))$.  We claim that
$C^*( \{ T_e, Q_v \} )$ is equal to $P C^*(d_r(E)) P$. Note that
if $P s_{\mu}s_{\nu}^*P = s_{\mu}s_{\nu}^* \ne 0 $ then
$s_{id}(\mu) = v_0$ and $s_{id}( \nu ) = w_0$ for some $v, w \in
E^0$ and $r_{id} ( \mu ) = r_{id} ( \nu )$.

If $r_{id} ( \mu ) \in E^0$, then $s_{\mu}s_{\nu}^*=
T_{\alpha}T_{\beta}^*$ for some paths $\alpha, \beta \in E^*$ and
hence $s_{\mu}s_{\nu}^* \in C^*(\{T_e,Q_v\})$. Suppose $r_{id}
(\mu) \notin E^0$.  Then $r_{id} (\mu) = r(e)_q$ for some $e \in
E^1$ with $1 \leq q \leq d_r( e )$ and there are $\alpha, \beta
\in E^*$  such that $s_\mu s_\nu^* = T_\alpha s_e
p_{f(r(e))_{d_r(e)}} s_e^* T_\beta^*$ if $q=d_r (e)$ and
\[
s_{\mu}s_{\nu}^*  = T_{\alpha} s_e s_{ f(r( e) )_{d_r ( e ) }}
\dots s_{f(r(e)_{q+1}} p_{f(r(e))_q} s_{f(r(e))_{q+1}}^* \dots
s_e^* T_{\beta}^*,
\]

\noindent otherwise. Since the vertices $f ( r ( e ) )_i$ for $2
\le i \le d_r ( r ( e ) )$ emit exactly one edge each, we have
$p_{f(r(e))_{i}} = s_{f(r(e))_{i-1}} s_{f(r(e))_{i-1}}^*$ and
hence we decrease $q$ in the expression for $s_\mu s_\nu^*$ until
we have $s_{\mu}s_{\nu}^* = T_{\alpha e}T_{\beta e}^* \in C^*(
T_e, Q_v )$ as required.

It remains to check that the corner is full.  To see this, we note
that $\Sigma H (E^0) = d_r(E)^0$ and apply Lemma~\ref{1.3}. Our
result follows.
\end{proof}

\begin{rmks}
\begin{itemize}
\item[(i)] Using in-delays we can convert row-finite graphs into
locally finite graphs (i.e. graphs where every vertex receives and
emits finitely many edges): If $E$ is row-finite and $v \in E^0$
receives edges $\{ e_i : i \in {\bf N} \}$, set $d_r ( v ) =
\infty$ and $d_r ( e_i ) = i$. If $v$ is a source we put $d_r (v)
= \infty$ and if $v$ receives finitely many edges we set $d_r
(v)=0$. Evidently $d_r : E^0 \cup E^1 \to {\bf N}$ is proper. The
resulting graph $d_r (E)$ is then locally finite with no sources.
Thus, combining Theorem \ref{outdelaysme} and Theorem
\ref{indelayme} we can show that for any graph $E$ there is a
locally finite graph with no sinks and sources $F$ such that $C^*
(E)$ is strongly Morita equivalent to $C^* (F)$.

\item[(ii)] An in-delay at a vertex $v$ with $d_r (v) \ge 1$
replaces $v \in E^0$  by the tree $\{ v_i : 0 \le i \le d_r (v) ,
f(v)_i : 1 \le i \le d_r (v) \}$ where $v$ is identified with the
leaf $v_0$. In combination with Remarks \ref{earlier} (iii) it
seems that we may get similar Morita equivalence results if we
replace vertices with more general trees (i.e.\ contractible
graphs) where the original vertex lies within the tree itself.

\item[(iii)] Not every in-delay can be expressed as an out-delay.
To see this observe that for the graph $E$ used in Examples
\ref{firstex} there can be no out-delay which corresponds to the
in-delay described in Example \ref{indelayexample}. It should not
be difficult to find examples where the graph contains no sources
and sinks.
\end{itemize}
\end{rmks}

\section{In-splittings} \label{is}

The following is adapted from \cite[Definition 2.4.7]{lm}: Let $E
= ( E^0 , E^1 , r , s )$ be a directed graph. For each $v \in E^0$
with $r^{-1}(v) \ne \emptyset$ partition the set $r^{-1} (v)$ into
disjoint nonempty subsets $\mathcal{E}^v_1 , \ldots ,
\mathcal{E}^v_{m(v)}$ where $m(v) \ge 1$ (if $v$ is a source then
we put $m(v)=0$). Let $\mathcal{P}$ denote the resulting partition
of $E^1$. We form the {\em in-split graph} $E_r({\mathcal P})$
from $E$ using the partition $\mathcal{P}$ as follows: Let
\begin{align*}
E_r ( \mathcal {P} )^0 &= \{ v_i : v \in E^0 , 1 \le i \le m(v) \}
\cup \{ v : m(v)  = 0 \} ,  \\
E_r ( \mathcal{P} )^1  &= \{ e_j : e \in E^1, 1 \le j \le m ( s (e) )
\} \cup \{ e : m (s(e)) = 0 \} ,
\end{align*}

\noindent
and define $r_{E_r ( \mathcal{P} )} , s_{E_r ( \mathcal{P} )} : E_r (
\mathcal{P}
)^1 \rightarrow E_r ( \mathcal{P} )^0$ by
\begin{align*}
s_{E_r ( \mathcal{P} )} ( e_j ) &= s(e)_j  \text{ and } s_{E_r (
\mathcal{P} )} ( e
) = s(e) \\
r_{E_r ( \mathcal{P} )} ( e_j ) &= r(e)_i \text{ and } r_{E_r (
\mathcal{P} )} ( e ) = r(e)_i \text{ where } e \in
\mathcal{E}^{r(e)}_i .
\end{align*}

\noindent Partition $\mathcal{P}$ is {\em proper} if for every
vertex $v$ which is a sink or emits infinitely many edges we have
$m(v)=0,1$. That is, we cannot in-split at a sink or vertex with
infinite valency.

To relate the graph algebras of a graph and its in-splittings we
use a variation of the method introduced in \cite[\S 4.2]{d}: If
$E_r ( \mathcal{P} )$ is the in-split graph formed from $E$ using
the partition $\mathcal{P}$ then we may define a Drinen
range-vector $d_{r , \mathcal{P}} : E^0 \cup E^1 \rightarrow {\bf
N} \cup \{ \infty \}$ by $d_{r , \mathcal{P}} ( v ) = m ( v )-1$
if $m(v) \ge 1$ and $d_{r, \mathcal{P}}(v)=0$ otherwise. For $e
\in \mathcal{E}^{r(e)}_i$ we put $d_{r , \mathcal{P}} ( e ) =
i-1$. Hence, if $v$ receives $n \ge 2$ edges then we create a
in-delayed graph in which $v$ is given delay of size $m(v)-1$ and
all edges with range $v$ are given a delay one less than their
label in the partition of $r^{-1} (v)$. If $v$ is a source or
receives only one edge then there is no delay attached to $v$.

\begin{examples} \label{isex}
\begin{itemize}
\item[(i)] Examples of proper in-splittings are found in
\cite[Figure 2.4.6]{lm} and \cite[\S 4.2]{d}.

\item[ii)] An in-splitting is not proper if we in-split at a sink,
such as for
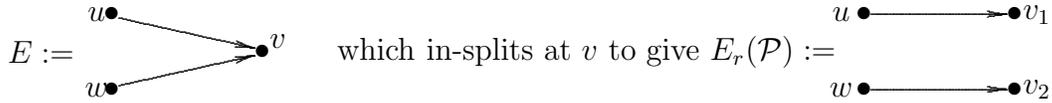
\[
\beginpicture

\setcoordinatesystem units <1cm,0.5cm>

\setplotarea x from -1.5 to 12, y from -1.5 to 1.5

\put{$E:=$}[r] at -1.5 0

\put{$\bullet$} at -1 1

\put{$\bullet$} at -1 -1

\put{$\bullet$} at 1 0

\arrow <0.25cm> [0.1,0.3] from -0.9 -0.9 to 0.9 -0.1

\arrow <0.25cm> [0.1,0.3] from -0.9 0.9 to 0.9 0.1

\put{which in-splits at $v$ to give $E_r ( \mathcal{P} ) :=$}[l]
at 2 0

\put{$\bullet$} at 9 -1

\put{$\bullet$} at 9 1

\put{$\bullet$} at 11 1

\put{$\bullet$} at 11 -1

\arrow <0.25cm> [0.1,0.3] from 9.1 1 to 10.9 1

\arrow <0.25cm> [0.1,0.3] from 9.1 -1 to 10.9 -1

\put{$u$} at -1.2 1

\put{$w$} at -1.2 -1

\put{$v$}[bl] at 1.1 0.1

\put{$u$} at 8.7 1

\put{$w$} at 8.7 -1

\put{$v_1$} at 11.3 1

\put{$v_2$} at 11.3 -1

\endpicture
\]

\noindent The associated in-delayed graph is
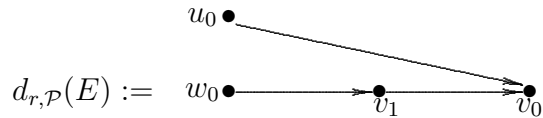
\[
\beginpicture

\setcoordinatesystem units <1cm,1cm>

\setplotarea x from -1 to 5, y from -0.2 to 1

\put{$d_{r , \mathcal{P}} ( E ) : =$}[r] at -1 0

\put{$\bullet$} at 0 1

\put{$\bullet$} at 0 0

\put{$\bullet$} at 2 0

\put{$\bullet$} at 4 0

\arrow <0.25cm> [0.1,0.3] from 0.1 0.9 to 3.9 0.1

\arrow <0.25cm> [0.1,0.3] from 0.1 0 to 1.9 0

\arrow <0.25cm> [0.1,0.3] from 2.1 0 to 3.9 0

\put{$u_0$} at -0.35 1

\put{$w_0$} at -0.35 0

\put{$v_1$} at 2.1 -0.2

\put{$v_0$} at 4 -0.2

\endpicture
\]
As $C^* ( E_r ( \mathcal{P} ) )$ has two ideals and $C^* (
d_{r,\mathcal{P}} (E) )$ one they are not Morita equivalent.

\item[(iii)] In-splittings at infinite valence vertices are not
proper, such as in
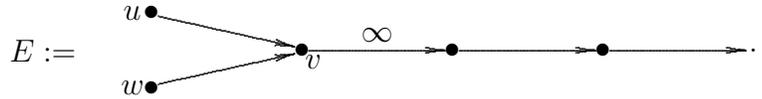
\[
\beginpicture

\setcoordinatesystem units <1cm,0.5cm>

\setplotarea x from -1 to 10, y from -0.5 to 1

\put{$E:=$}[r] at -1 0

\put{$\bullet$} at 0 1

\put{$\bullet$} at 0 -1

\put{$\bullet$} at 2 0

\put{$\bullet$} at 4 0

\put{$\bullet$} at 6 0

\put{$\ldots$} at 8 0

\arrow <0.25cm> [0.1,0.3] from 0.1 0.9 to 1.9 0.1

\arrow <0.25cm> [0.1,0.3] from 0.1 -0.9 to 1.9 -0.1

\arrow <0.25cm> [0.1,0.3] from 2.1 0 to 3.9 0

\arrow <0.25cm> [0.1,0.3] from 4.1 0 to 5.9 0

\arrow <0.25cm> [0.1,0.3] from 6.1 0 to 7.9 0

\put{$\infty$} at 3 0.4

\put{$u$} at -0.25 1

\put{$w$} at -0.25 -1

\put{$v$} at 2.15 -0.3

\endpicture
\]

\noindent which in-splits at $v$ to give $E_r ( \mathcal{P} ) :=$
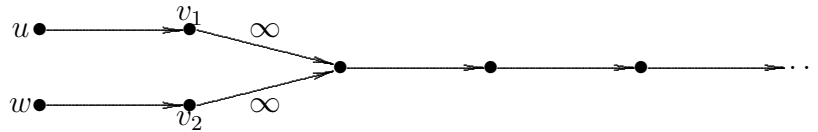
\[
\beginpicture

\setcoordinatesystem units <1cm,0.5cm>

\setplotarea x from -1 to 10, y from -1.4 to 1.4

\put{$\bullet$} at 0 -1

\put{$\bullet$} at 0 1

\put{$\bullet$} at 2 1

\put{$\bullet$} at 2 -1

\put{$\bullet$} at 4 0

\put{$\bullet$} at 6 0

\put{$\bullet$} at 8 0

\put{$\ldots$} at 10.2 0

\arrow <0.25cm> [0.1,0.3] from 0.1 1 to 1.9 1

\arrow <0.25cm> [0.1,0.3] from 0.1 -1 to 1.9 -1

\arrow <0.25cm> [0.1,0.3] from 2.1 -1 to 3.9 -0.1

\arrow <0.25cm> [0.1,0.3] from 2.1 1 to 3.9 0.1

\arrow <0.25cm> [0.1,0.3] from 4.1 0 to 5.9 0

\arrow <0.25cm> [0.1,0.3] from 6.1 0 to 7.9 0

\arrow <0.25cm> [0.1,0.3] from 8.1 0 to 9.9 0

\put{$u$} at -0.25 1

\put{$w$} at -0.25 -1

\put{$v_1$} at 2 1.4

\put{$v_2$} at 2 -1.4

\put{$\infty$} at 3 1

\put{$\infty$} at 3 -1

\endpicture
\]

\noindent The associated in-delayed graph is
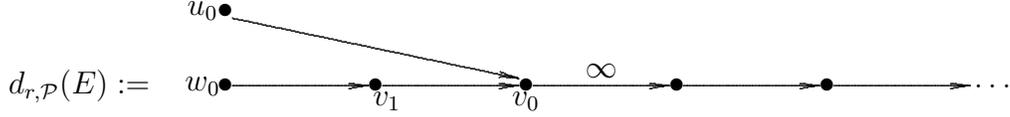
\[
\beginpicture

\setcoordinatesystem units <1cm,1cm>

\setplotarea x from -1 to 10, y from -1 to 1

\put{$d_{r , \mathcal{P}} ( E ) : =$}[r] at -1 0

\put{$\bullet$} at 0 0

\put{$\bullet$} at 0 1

\put{$\bullet$} at 2 0

\put{$\bullet$} at 4 0

\put{$\bullet$} at 6 0

\put{$\bullet$} at 8 0

\put{$\ldots$} at 10.2 0

\arrow <0.25cm> [0.1,0.3] from 0.1 0.9 to 3.9 0.1

\arrow <0.25cm> [0.1,0.3] from 0.1 0 to 1.9 0

\arrow <0.25cm> [0.1,0.3] from 2.1 0 to 3.9 0

\arrow <0.25cm> [0.1,0.3] from 4.1 0 to 5.9 0

\arrow <0.25cm> [0.1,0.3] from 6.1 0 to 7.9 0

\arrow <0.25cm> [0.1,0.3] from 8.1 0 to 9.9 0

\put{$u_0$} at -0.3 1

\put{$w_0$} at -0.3 0

\put{$v_1$} at 2.15 -0.2

\put{$v_0$} at 4 -0.2

\put{$\infty$} at 5 0.2

\endpicture
\]
In this case $C^* ( E_r ( \mathcal{P} ) )$ has two ideals, whereas
$C^*(d_r (E))$ only has one. Thus these algebras are not Morita
equivalent.
\end{itemize}
\end{examples}

\begin{rmk} \label{properpoint}
If $\mathcal{P}$ is proper then every vertex $v$ which is either a
sink or a vertex of infinite valency occurs only as $v$ or $v_1$
in $E_r ( \mathcal{P} )$ and only as $v_0$ in $d_{r ,\mathcal{P}}
(E)$. In particular, if $\mathcal{P}$ is a proper partition and
$v$ is a sink or infinite valence vertex, then there no edges of
the form $e_j$ for $j \ge 2$ with $s(e)=v$ in $E_r ( \mathcal{P}
)$ and no edges of the form $f(v)_i$ in $d_{r , \mathcal{P}} (E)$.
\end{rmk}

\begin{thm} \label{insplitindelay}
Let $E$ be a directed graph, $\mathcal{P}$ a partition of $E^1$,
$E_r ( \mathcal{P} )$ the in-split graph formed from $E$ using
$\mathcal{P}$ and $d_{r , \mathcal{P}} : E^0 \cup E^1 \rightarrow
{\bf N} \cup \{ \infty \}$  the Drinen range-vector defined as
above. Then $C^* ( E_r ( \mathcal{P} ) ) \cong C^* ( d_{r
,\mathcal{P}} (E) )$ if and only if $\mathcal{P}$ is proper.
\end{thm}

\begin{proof}
Let $\{ s_f , p_w : f \in d_{r , \mathcal{P}} (E)^1 ,  w \in d_{r
, \mathcal{P}} (E)^0  \}$ be a Cuntz-Krieger $d_{r, \mathcal{P} }
(E)$-family. To simplify our definitions, for $v \in E^0$ we put
$s_{f(v)_{0}} = p_{v_{0}}$. For $e \in E_r ( \mathcal{P} )^1$ we
define $T_e = s_e$. For $e_j \in E_r ( \mathcal{P} )^1$ with $1
\le j \le m (s(e))$ we define
\[
T_{e_j}  = s_{f(s(e))_{j-1}} \ldots s_{f(s(e))_1} s_e  .
\]

\noindent For $1 \le i \le m ( v )$, define $Q_{v_i} =
p_{v_{i-1}}$; if $m(v)=0$, define $Q_v = p_{v_0}$. Then $\{ T_g ,
Q_u : g \in E_r ( \mathcal{P} )^1 , u \in E_r ( \mathcal{P} )^0
\}$ is a Cuntz-Krieger $E_r ( \mathcal{P} )$-family with $Q_u \neq
0$ for all $u$.

Let $\{ t_g , q_u \}$ be the canonical generators of $C^* ( E_r (
\mathcal{P} ) )$. By the universal property of $C^* ( E_r (
\mathcal{P} ))$ there is a homomorphism $\pi : C^* ( E_r (
\mathcal{P} ) ) \to C^* ( d_{r , \mathcal{P}} (E) )$ such that
$\pi ( t_g ) = T_g$ and $\pi ( q_u ) = Q_u$. We claim that $\pi$
is surjective, that is $\{ T_g, Q_u \}$ generates $C^*( d_{r ,
\mathcal{P} } ( E ))$.

For $w \in d_{ r  \mathcal{P}} (E )^0$ we have $p_w \in C^*( \{
T_g , Q_u \} )$ by definition. For $e \in E_r ( \mathcal{P} )^1$
we have $s_e = T_e \in C^* ( \{ T_g , Q_u \} )$. Since
$\mathcal{P}$ is proper, by Remark \ref{properpoint} there are no
edges in $d_{r , \mathcal{P}} (E)$ of the form $f(r(e))_j$ with
$r(e)$ a sink. In particular, every edge $f(v)_j$ in $d_{r ,
\mathcal{P}} (E)$ is of the form $f(s(e))_j$ where $s(e)$ has
finite valency. For $1 \le j \le m(s(e)) - 1 = d_{r , \mathcal{P}}
(s(e))$
\[
T_{e_{j+1}} T_{e_j}^* = s_{f(s(e))_{j}} \dots s_{f(s(e))_{1}} s_e
s_e^* s_{f(s (e))_{1}} \dots s_{f(s (e))_{j-1}}  \\
\]

\noindent and since $v= s(e)$ has finite valency we have
\begin{align*}
\sum_{s(e) = v} T_{e_{j+1}} T_{e_j}^*  &= s_{f(s(e))_{j}} \dots
s_{f(s(e))_{1}}  \left(\sum_{s(e) = v_0} s_e s_e^*\right)
s_{f(s(e))_{1}}^* \dots s_{f(s(e))_{j-1}}^*  \\
&= s_{f(s(e))_{j}} \dots s_{f(s(e))_1} p_{v_0} s_{f(s(e))_1}^*
\dots s_{f(s(e))_{j-1}}^* = s_{f(s(e))_{j}}.
\end{align*}

\noindent Then $s_{f(s(e))_j} \in C^*( T_g ,Q_u  )$ and our claim
follows.

For $z \in \bf{T}$ define an action $\alpha$ on $C^*(d_{r,\mathcal
P}(E))$ by $\alpha_z ( p_v ) = p_v$ for $v \in d_{r ,{\mathcal P}}
(E)^0$,  $\alpha_z (s_{e} ) = z s_e$ for $e \in E_r ( \mathcal{P}
)^1$, and $\alpha_z (s_{f(v)_i} ) = s_{f(v)_i}$ for $1 \le i \le
d_{r, \mathcal{P}} (v)$. Since $\gamma \circ \pi = \pi \circ
\alpha$ where $\gamma$ is the usual gauge action on
$C^*(E_r({\mathcal P}))$, by Theorem \ref{bhrsz}
$C^*(d_{r,\mathcal P}(E)) \cong C^* ( E_r ( \mathcal{P} ) )$.

If $\mathcal{P}$ is not proper then there is a non-trivial
in-splitting at a sink or a vertex of infinite valency. The
 graph $E_r ( \mathcal{P} )$ will have at least one
more sink or vertex of infinite valency than $d_{r , \mathcal{P}}
(E)$ and hence $C^* ( E_r ( \mathcal{P} ) )$ will have more ideals
than $C^* ( d_{r , \mathcal{P}} (E) )$.
\end{proof}

\noindent Applying Theorem \ref{insplitindelay} and Theorem
\ref{indelayme} we have:

\begin{cor} \label{insplitme}
Let $E$ be a directed graph, $\mathcal{P}$ a partition of $E^1$
and $E_r ( \mathcal{P} )$ the in-split graph formed from $E$ using
$\mathcal{P}$, then $C^* ( E_r ( \mathcal{P} ) )$ is strongly
Morita equivalent to $C^* ( E )$ if and only if $\mathcal{P}$ is
proper.
\end{cor}

\section{Connections with Strong Shift Equivalence} \label{shifteq}
In \cite{b1} the following definition (which generalises one given in \cite{ash}) was given for
{\em elementary strong shift equivalence} of directed graphs which contain no sinks.

\begin{dfn} \label{strongshift}
Let $E_i = (E^0_i, E^1_i,r^i,s^i)$ for $i = 1,2$ be directed
graphs.  Suppose there is a directed graph
$E_3=(E^0_3,E^1_3,r_3,s_3)$ such that:
\begin{itemize}
\item[a)] $E^0_3 = E^0_1 \cup E^0_2$ and $E^0_1 \cap E^0_2 =
\emptyset$. \item[b)] $E^1_3 = E^1_{12} \cup E^1_{21}$ where
$E^1_{ij} := \{ e \in E^1_3 : s_3(e)\in E^0_i, r_3 (e) \in
E^0_{j}\}$. \item[c)] For $i= 1,2$ there are range and
source-preserving bijections $\theta_i : E^1_i \to
E^2_3(E^0_i,E^0_i)$ where for $i \in \{1,2\}$, $E^2_3(E^0_i,E^0_i)
:= \{ \alpha \in E^2_3 : s_3(\alpha) \in E^0_i, r_3(\alpha) \in
E^0_i\}$.
\end{itemize}
Then we say that $E_1$ and $E_2$ are {\em elementary strong shift equivalent}
$(E_1 \thicksim_{ES} E_2)$ via $E_3$.
\end{dfn}

\noindent The equivalence relation $\thicksim_S$ on directed
graphs generated by elementary strong shift equivalence is called
\emph{strong shift equivalence}. Row-finite graphs which are
strong shift equivalent have Morita equivalent $C^*$-algebras (see
\cite[Theorem 5.2]{b1}).

\begin{prop} \label{outsplitesse} Let $E$ be a directed graph with no sinks
and $E_s(\mathcal{P})$ be an
out-split graph formed from $E$ using ${\mathcal P}$.  Then $E
\thicksim_{ES} E_s({\mathcal P})$.
\end{prop}

\begin{proof}
One constructs a bipartite graph $E_3$ in the following manner.
Let $E^0_3 = E^0 \cup E_s({\mathcal P})^0$. For each $v \in E^0$
draw an edge $e^i_v$ to the corresponding split vertices $v^i \in
E_s({\mathcal P})^0$ with $s(e^i_v) = v$ and $r(e^i_v) = v^i$. For
each set of edges $\{e^i\}_{i=1}^{m(r(e))} \subseteq E_s(P)^1$
with $s(e^i) = v^i$ and $r(e^i) \in \{w^j\}_{j=1}^{m(w)}$, draw an
edge $e^i_{v,w}$ with $s(e^i_{v,w}) = v^i$ and $r(e^i_{v,w}) = w$.
The graph $E_3$ satisfies the conditions of
Definition~\ref{strongshift} and hence $E \thicksim_{ES}
E_s({\mathcal P})$ via $E_3$.
\end{proof}

\noindent In a similar manner we may show:

\begin{prop} \label{insplitesse} Let $E$ be a directed graph with
no sinks and $E_r(\mathcal{P})$ be an in-split
graph formed from $E$ using ${\mathcal P}$.  Then
$E \thicksim_{ES} E_r({\mathcal P})$.
\end{prop}

\begin{rmk}
Proposition~\ref{insplitesse} and \cite[Theorem 5.2]{b1} enable us
to give another proof that the $C^*$-algebras of a row-finite
directed graph and its in-splitting are Morita equivalent. We have
analogous results for in-amalgamations and out-amalgamations as
they are the reverse operations of in-splittings and
out-splittings.
\end{rmk}


\begin{thebibliography}{BPRSz}

\bibitem[Ar]{ar} R.J.~Archbold.
\newblock On the `Flip-Flop' Automorphism of $C^* ( S_1 , S_2 )$.
\newblock {\em Quart.\ J.\ Math} {\bf 30}: 129--132 (1979).

\bibitem[Ash]{ash} B.~Ashton.
\newblock Morita equivalence of graph $C^*$-algebras.
\newblock Honours Thesis.
\newblock Uni.\ of Newcastle, 1996.

\bibitem[B1]{b1} T.~Bates.
\newblock Applications of the gauge--invariant uniqueness theorem for
the Cuntz--Krieger algebras of directed graphs.
\newblock Bull. Austral. Math. Soc., {\bf 65}: 57--67 (2002).

\bibitem[B2]{ba2} T.~Bates,
\newblock On the primitive ideal spaces of row-finite graphs.
\newblock {\em Preprint: University of New South Wales, 2001.}

\bibitem[BHRSz]{bhrsz} T.~Bates, D.~Pask, I.~Raeburn, W.~Szyma\'{n}ski.
\newblock The ideal structure of the $C^*$--algebras of infinite graphs.
\newblock {\em Illinois J.\ Math.}, {\bf 46}: 1159--1176, 2002.


\bibitem[BPRSz]{bprs} T.~Bates, D.~Pask, I.~Raeburn, W.~Szyma\'{n}ski.
\newblock The $C^*$--algebras of row--finite graphs.
\newblock {\em New York J.\ Math.}, {\bf 6}: 307--324, 2000.

\bibitem[Br]{bb} B.~Brenken.
\newblock $C^*$--algebras of infinite graphs and Cuntz--Krieger algebras.
\newblock {\em Can.\ Math.\ Bull.}, to appear

\bibitem[CG]{cg} T.~Crisp and D.~Gow.
\newblock Contractible subgraphs and Morita equivalence of graph $C^*$-algebras,
{\em Preprint: University of Newcastle, 2003}.

\bibitem[CK]{ck} J.~Cuntz and W.~Krieger.
\newblock A class of $C^*$-algebras and topological Markov chains.
\newblock {\em Invent. Math.} {\bf 56}: 251--268, 1980.

\bibitem[D]{d} D.~Drinen.
\newblock Flow equivalence and graph groupoid isomorphism.
\newblock {\em Preprint: Dartmouth College, 2001.}

\bibitem[DS]{ds} D.~Drinen and N. Sieben.
\newblock $C^*$-equivalences of graphs.
\newblock {\em J. Operator Theory} {\bf 45}: 209--229 (2001).

\bibitem[DT]{dt} D.~Drinen and M.~Tomforde.
\newblock The $C^*$-algebras of arbitrary graphs.
\newblock {\em Rocky Mountain J.\  Math.}, to appear

\bibitem[EW]{ew} M.~Enomoto,  Y.~Watatani.
\newblock  A graph theory for $C^*$--algebras.
\newblock {\em Math. Japon.} {\bf 25}: 435--442, 1980.

\bibitem[FW]{fw} M.~Fujii,  Y.~Watatani.
\newblock  Cuntz--Krieger algebras associated with adjoint graphs.
\newblock {\em Math. Japon.} {\bf 25}: 501--506, 1980.

\bibitem[GT]{gt} J.L.~Gross and T.W.~Tucker.
\newblock {\em Topological graph theory.}
\newblock Wiley Interscience Series in Discrete Mathematics
and Optimization, First edition, (1987).

\bibitem[KQR]{kqr} S.~Kaliszewski, J.~Quigg, and I.~Raeburn.
\newblock Skew products and crossed products by coactions.
\newblock {\em J. Operator Theory} {\bf 46}: 411-433, (2001).

\bibitem[KPRR]{kprr} A.~Kumjian, D.~Pask, I.~Raeburn, and J.~Renault.
\newblock Graphs, groupoids and Cuntz--Krieger algebras.
\newblock{\em J. Funct. Anal.} {\bf 144}: 505--541, 1997.

\bibitem[KPR]{kpr} A.~Kumjian, D.~Pask, I.~Raeburn.
\newblock Cuntz--Krieger algebras of directed graphs,
\newblock {\em Pacific. J. Math.}, {\bf 184}: 161--174, (1998).

\bibitem[KP]{kp} A.~Kumjian, D.~Pask.
\newblock $C^*$-algebras of directed graphs and group actions,
\newblock {\em Ergod. Th. \& Dynam. Sys.}, {\bf 19}: 1503--1519, (1999).

\bibitem[LM]{lm} D.~Lind and B.~Marcus.
\newblock {\em An introduction to symbolic dynamics and coding.}
\newblock Cambridge University Press., 1995.

\bibitem[MRS]{mrs} M.H.~Mann, I.~Raeburn, C.E.~Sutherland.
\newblock Representations of finite groups and Cuntz--Krieger
algebras
\newblock {\em Bull. Austral. Math. Soc.} {\bf 46}: 225--241,
1992.

\bibitem[PR]{pr} D.~Pask, and I.~Raeburn.
\newblock On the $K$-theory of Cuntz--Krieger algebras
\newblock {\em Proc. RIMS Kyoto} {\bf 32}: 415--443, 1996.

\bibitem[PS]{ps} W.~Parry, and D.~Sullivan.
\newblock A topological invariant of flows on 1-dimensional spaces.
\newblock {\em Topology} {\bf 14}: 297--299, 1975.

\bibitem[RS]{rs} I.~Raeburn and W~Szyma\'{n}ski.
\newblock Cuntz--Krieger algebras of infinite graphs and matrices.
\newblock {\em Trans. \ Amer.\ Math.\ Soc.}, to appear.

\bibitem[St]{st} J.~Stillwell.
\newblock {\em Classical Topology and combinatorial group theory.}
\newblock Volume~{\bf 72} of {\em Graduate Texts in
Mathematics}, Springer--Verlag, (1980).

\bibitem[T]{t} M.~Tomforde.
\newblock A unified approach to Exel-Laca algebras and $C^*$-algebras associated
to graphs.
\newblock {\em J.\ Operator Theory}, to appear.

\bibitem[W]{wi} R.F.~Williams.
\newblock Classification of subshifts of finite type.
\newblock Ann. of Math. (2) {\bf 98}: 120--153, 1973.
\end{thebibliography}
\end{document}